\documentclass[a4paper, 11pt]{amsart}

\usepackage[T1]{fontenc}
\usepackage[utf8]{inputenc}
\usepackage{geometry}
\usepackage{amssymb}

\makeatletter
\@namedef{subjclassname@2020}{%
  \textup{2020} Mathematics Subject Classification}
\makeatother

\allowdisplaybreaks[2]

\newtheorem{theorem}{Theorem}[section]
\newtheorem{proposition}[theorem]{Proposition}
\newtheorem{lemma}[theorem]{Lemma}

\theoremstyle{definition}

\newtheorem{remarks}[theorem]{Remarks}
\newtheorem{definition}[theorem]{Definition}

\numberwithin{equation}{section}

\newcommand\NN{\mathbb{N}}
\newcommand\CC{\mathbb{C}}
\newcommand\ZZ{\mathbb{Z}}
\begin{document}
\title[ Waring problems across algebra]{ Waring problems across algebra}

\author{Matej Brešar}
\address{Faculty of Mathematics and Physics, University of Ljubljana \&
Faculty of Natural Sciences and Mathematics, University of Maribor \& IMFM, Ljubljana, Slovenia}
\email{matej.bresar@fmf.uni-lj.si}

\author{Consuelo Mart\' inez}
\address{Departamento de Matem\' aticas, Universidad de Oviedo,
Spain}
\email{cmartinez@uniovi.es}
   \thanks{The first author was partially supported by the ARIS Grants P1-0288 and J1-60025.  
The second author was partially supported bt the Spanish project MCIU-22-PID2021-123461NB-C22.}
\date{}

\begin{abstract} 
The paper surveys various Waring type problems in groups, Lie algebras, and associative algebras.  
\end{abstract}

\maketitle

\section{Introduction}

The classical Waring problem, proposed by Edward Waring in 1770, asks whether  for every
positive integer $n$ there exists a positive integer $N$ such that
every positive integer $b$ can be expressed as
$$ b = a_1^n + \cdots + a_N^n$$
for some nonnegative integers $a_i$. In 1909, David Hilbert
gave a positive answer to this question.

Denote the smallest such
number $N$ by $g(n)$. 
 Lagrange's four-square theorem states that
 $g(2)=4$. The value of the function $g$ is also known for  many other  $n$, for example,  $g(3)=9$,
 $g(4)=19$, $g(5)=37$, etc. However, not everything concerning the function $g$ is understood, and there are also 
  other problems arising from the classical Waring problem that are still an active research area in Number Theory. 

The problems and results outlined above
inspired analogous problems in various algebraic structures. The classical situation has been extended in many different directions, in particular by replacing the role of the powers $a^n$ by more general expressions.

The purpose of this paper is to survey several Waring type results for groups, Lie algebras, and associative algebras. It should be mentioned that these research areas are too broad to be covered in a short paper. This especially holds for the group-theoretic part of this line of research (see, for example,  \cite{Se} and \cite{GKP}). 
We will  focus on topics that are close to our research interests.

\bigskip

\section{Waring problems in groups}

Waring type problems in group theory are linked to the notion of  {\it width} of a word $w$. A word $w = w(x_1, \ldots, x_n)$ is a nontrivial element ($w \neq 1$) in the free group $F(\infty)$ of countable rank (we can also see the word as an element of the free group of rank $n$, $F_n$).

Given a group $G$ we can consider the set of all elements of $G$ obtained by substituting arbitrary elements of $G$ in the word $w$: $$w(G) = \{w(g_1, \ldots, g_n) \,|\, g_1, \ldots, g_n \in G\}.$$

The {\it verbal set} $w(G)$ generates the {\it verbal subgroup} $\langle w(G) \rangle$.
For instance, if $w =[x,y] = x^{-1}y^{-1}xy$, then $\langle w(G) \rangle = G'$, the derived subgroup of $G$.

\begin{definition}
The word $w$ is {\it elliptic} on $G$ if $\langle w(G) \rangle = \underbrace{ w(G)^{\pm 1} \cdots w(G)^{\pm 1}}_N$ for some $N$.  
The smallest $N$ satisfying the above condition is called the ${\it width}$ of $w$ in $G$.  
\end{definition}
See \cite{Se} for a complete survey on elliptic words in groups.

\medskip

\subsection{Finite groups}

The above definition holds for all groups. It is clear that every word has finite width in a finite group.  So, why should we be interested in this notion for finite groups?  There are many interesting problems that can be considered.  We may ask different questions about distributions of words in finite groups, or we may ask about the existence of uniform bounds  for every finite group or for all groups in some particular class of finite groups.

Notice that if $G$ is a simple group then either $ w(G) = 1$ or $\langle w(G) \rangle = G$, since the verbal subgroup is always a normal subgroup.

Ore considered the word $w =[x,y]$ and proved that it has width 1 in every alternating group $A_m$, $m \geq 5$, and conjectured that the same  is true for any simple finite group.  This conjecture, known as {\it Ore conjecture}, inspired a lot of work.

Saxl and Wilson proved in \cite{SW} 
(see also \cite{Wilson}) that there is a bound for the width of the commutator in any simple finite group. This result was extended in 2001 by Liebeck and Shalev proving that  for an arbitrary word $w$ there is a bound, $c = c(w)$, (that depends only of $w$) for the width of $w$ in an arbitrary simple finite group.  However, no explicit value for this bound is given (see \cite{LSh}).

In 2009, Shalev \cite{Shalev} proved that if the order of a simple finite group $G$ is big enough, the width of an arbitrary word $w$ in $G$ is bounded by 3. 

In 2010, Liebeck, O'Brien, Shalev and Tiep proved the Ore conjecture:  Every element in a simple finite group is a commutator \cite{OreCon}.  The proof uses character theory, algebraic geometry and computational methods.

Further, Larsen, Shalev and Tiep \cite{LST} proved that, for an arbitrary word $w$, any element in a finite non-abelian simple group of sufficiently high order can be written as the product of two elements in $w(G)$.

\subsection{Verbal ellipticity of groups}

We have already mentioned that every word has finite width in a finite group.  Thus, every finite group $G$ is {\it verbally elliptic}, i.e.,
every word is elliptic on $G$.  Are there infinite {\it verbally elliptic} groups?
There are several results with regard to this question.

Romankov proved in \cite{Ro} that the following groups are {\it verbally elliptic}:
\begin{enumerate}
    \item[1.] Finitely generated virtually nilpotent groups.
    \item[2.] Virtually abelian groups.
    \item[3.] Finite rank virtually nilpotent groups.
\end{enumerate}
He also constructed an example of a solvable linear group on which the word $w = [[x_1,x_2],[x_3,x_4]]$ is not elliptic.

N. Nikolov and D. Segal (\cite{NS}, \cite{NS1}, \cite{Se}) proved that there exists a function $f(m,n)$ such that the word $w = x^n$ has width $\leq f(m,n)$ on every finite $m$-generated group.

This result, together with the solution of the Restricted Burnside Problem implies the proof of the conjecture of J.\,P. Serre: {\it a subgroup of finite index in a finitely generated profinite group is closed}. 

A. Jaikin Zapirain proved in \cite{J-Z} that a finitely generated $p$-adic analytic group is verbally elliptic.

The Nottingham group over a finite field is also verbally elliptic, as proved by B. Klopsch in \cite{Kl}.  The same  was proved for the Nottingham group over a field of zero characteristic in \cite{MM}.

\medskip

\subsection {pro-$p$ groups}

In this section we will consider {\it pro-$p$ groups}. Let us first recall some notions from group theory (see \cite{Se}).

A group $G$ is a {\it residually p-group} if
$$ \cap\{H \triangleleft G \;| \;|G:H| = p^k, k \geq 1 \} = (1).$$
The above subgroups  $H$ can be considered as a neighborhood basis of the element $1$, making $G$ a topological group.
If the corresponding topology is complete, then $G$ is a {\it pro-$p$ group}.  Otherwise, the completion $G_p$ of the group $G$ is a pro-$p$ group that is called the {\it pro-$p$ completion} of $G$.

Let $F(\infty)$ be, as before, the free group on countable many generators $x_1, x_2, x_3, \ldots$  Its  pro-$p$ completion, $F = F(\infty)_p$, is the free  pro-$p$ group on $x_1, x_2, x_3, \ldots$.  An element $ w \in F$ is a {\it word} if it involves finitely many generators.  

Let $G$ be a  pro-$p$ group and $w \in F$ a word.  Let us consider again the set $w(G) = \{w(g_1, \ldots, g_n) | g_1, \ldots, g_n \in G\}$ and denote by $\overline{ \langle w(G) \rangle}$ the closed subgroup generated by $w(G)$.

We say that the word $w \in F$ is {\it elliptic} on $G$ (or has finite {\it verbal width} on $G$) if there is an integer $N \geq 1$ such that $\overline{ \langle w(G) \rangle} = \underbrace{w(G)^{\pm 1} \cdots w(G)^{\pm 1}}_N$.
This is equivalent to say that the discrete subgroup of $G$ generated by $w(G)$ is closed.

In 2018, the following result was proved in \cite{M}.

\begin{theorem} Let $\Gamma$ be a finitely generated residually-p torsion group and $G = \Gamma_p$ its pro-$p$-completion. Then an arbitrary word $w \in F$ is elliptic on $G$.
\end{theorem}

Note that important groups as Golod-Shafarevich, Grigorchuk or Gupta-Sidki groups are residually-$p$ torsion groups.

In the proof of this result the following important result by E. Zelmanov (see \cite {Z}) was used.

\begin{theorem}  Let $G$ be a  pro-$p$ group satisfying a nontrivial pro-$p$-identity.  If $G$ has a dense finitely generated torsion discrete subgroup, then $G$ is finite.
\end{theorem}

Given a  pro-$p$ group $G$ and $1 \neq w \in F$, the group $G$ satisfies the pro-$p$ identity $w = 1$ if $w(g) = 1$.

In the same paper \cite{M} the notion of strong ellipticity is considered. Let us consider a word $w = w(x_1, \cdots, x_n) \in F$ and $G$ a  pro-$p$ group. Choose $n-1$ elements  $a_1,\ldots, a_{i-1}, a_{i+1},\ldots, a_n \in G$ and fix all variables $x_j = a_j$, except $x_i$.  We construct the set $w(G,i,\alpha)$, where $\alpha = (a_1,\ldots, a_{i-1},a_{i+1},\ldots,a_n)$, by
$$w(G,i,, \alpha) = \{w(a_1, \ldots, a_{i-1},g,a_{i+1}, \ldots, a_n) \; | \; g \in G\}.$$

\begin{definition}  A word $w \in F$ is strongly elliptic on $G$ if there exist finite subsets $M_i \subset \underbrace{ G \times \cdots \times G}_{n-1}$, $1 \leq i \leq n$ and an order on the set $\cup_{i=1}^n M_i = \{\alpha_1 < \alpha_2 < \cdots <\alpha_q\}$, $\alpha_k \in M_{i_k}$, $1 \leq i_k \leq n$, such that the verbal subgroup $\overline{\langle w(G) \rangle}$ is equal to
$$ \overline{\langle w(G) \rangle} = w(G,i_1,\alpha_1)^{\pm 1} \cdots w(G,i_q,\alpha_q)^{\pm 1}.$$
\end{definition}

It is clear that any strongly elliptic word on $G$ is also elliptic.

In order to prove that the derived subgroup of a finitely generated  pro-$p$ group is closed, P. Serre proved (see \cite{Serre})  that the commutator is strongly elliptic on finitely generated  pro-$p$ groups.

The mentioned result in \cite{M} can be extended to {\it multilinear words}.  Let us define this notion.

As before let
 $F$ be the free  pro-$p$ group. Consider its central series,  $F = F_1 \geq F_2 \geq \cdots$. 

\begin{definition} Let  $w \in F$ be a word.  Supose that $w \in F_n \setminus F_{n+1}$.  The word $w$ is {\it multilinear} if $w = \tilde w w'$, where

(i)  $\tilde w$ is a product of left-normed commutators
$$\tilde w = \prod_{\sigma \in S_n} [ \cdots[x_{\sigma(1)}, x_{\sigma(2)}], \cdots, x_{\sigma(n)}]^{k_{\sigma}},$$
where  $k_{\sigma} \in Z_p$ are $p$-adic intergers such that the element $\tilde w F_{n+1}$ does not belong to $pL(F)$, and

(ii) the element $w' \in F_{n+1}$ is a converging product of commutators in $x_1, \ldots, x_n$ of length $\geq n+1$ and each commutator involves all $n$ generators $x_1, \ldots, x_n$.
\end{definition}

Using the ideas of the well known linearization process in algebras, the following result can be proved:  For an arbitrary word $w \in F$ there exists a multilinear word $\overline x$ that lies in the closed subgroup generated by $w(F)$  (see \cite{M}).

The following theorem is also proved in \cite{M}.

\begin{theorem}   Let $\Gamma$ be a finitely generated residually-p torsion group and $G = \Gamma_p$ its pro-$p$-completion. Then an arbitrary multilinear word $w \in F$ is strongly elliptic on $G$.
\end{theorem}

\bigskip

\section{Waring problems in Lie algebras}

In many cases, Waring problems in Lie algebras are connected to similar problems in groups related to these Lie algebras. We will center our attention in the results of this line of work, which, however, is not the only one (see \cite{BGKP} or \cite{KHR}). For instance, in \cite{KHR} authors consider the bracket width in current Lie algebras. Let us briefly survey their results.

If $L$ is a Lie algebra over an infinite field $K$, $z \in [L,L]$, $l(z)$ is the smallest number $n$ such that $z$ can be expressed as a sum of $n$ Lie brackets: $z = \sum_{i=1}^n [x_i,y_i]$. The {\it bracket width} of $L$ is the supremum of lengths $l(z)$, where $z$ runs over the derived algebra $[L,L]$.

If $K$ is an algebraically closed field of characteristic zero, $L$ is a finite dimensional simple Lie $K$-algebra and $A$ is a commutative associative unital $K$-algebra, then the {\it current algebra} corresponding to $L$ and $A$ is the Lie algebra tensor product $L \otimes_K A$ with the product
$$[x \otimes a, y \otimes b] =[x,y] \otimes ab.$$

The authors proved the following results:

\begin{theorem} The {\it bracket width} of $L \otimes_K A$ is less than or equal to 2.
\end{theorem}

In the particular case that $A = K[[t]]$,  the result result obtained is:

\begin{theorem} 1)  If $L = {\rm sl}_2$, then the bracket width of $L \otimes_K K[[t]]$ is equal to 1.

2)  The bracket width of $L\otimes_K K[[t]]$ is equal to 2 if $L$ is of type $A_n$ or $C_n$ ($n \geq 2$).
\end{theorem}

\medskip

\subsection{Lie algebras linked to  pro-$p$ groups}

The concepts of ellipticity and strong ellipticity of groups can be defined for Lie algebras.

To start we recall the Lie rings linked to a  pro-$p$ group $G$. 

Consider the lower central series of $G$: $G = G_1 \geq G_2 \geq \cdots$, where $G_{i+1} = [G, G_i]$.
Each factor $G_i /G_{i+1}$ is abelian and can be considered as a module over $Z_p$, the ring of $p$-adic integers.
The direct sum
$$ L(G) = \sum_{i \geq 1} G_i/G_{i+1}$$
is a Lie algebra over the ring $Z_p$, with the Lie product being the linear extension to arbitrary elements of $L(G)$ of the product defined on homogeneous elements by
$$[aG_{i+1}, bG_{j+1}] = [a,b]G_{i+j+1}, \; a \in G_i, \; b \in G_j.$$

In the particular case of $G = F$, the free  pro-$p$ group, the Lie algebra $L(F)$ is the free Lie $Z_p$-algebra over the set of free generators $X = \{ {\bar x_i} = x_iF_2 \, | \, i \geq 1\}$.

Now we can give the corresponding definitions.

Let $Lie \langle X \rangle$ be the free Lie $Z_p$-algebra over the set of free generators $X = \{x_i \, |\, i\geq 1\}$.
Let $f(x_1, \ldots, x_n)$ be a multilinear element of $Lie \langle X \rangle$.  For an arbitrary Lie $Z_p$-algebra $L$ we can consider $f(L) = \{f(a_1, \ldots, a_n) \; | \; a_i \in L \}$.

\begin{definition}  The element $f$ is {\it elliptic} on $L$ if there is a natural number $N \geq 1$ such that the $Z_p$-linear span of $f(L)$ is equal to
$$Span(f(L)) = \underbrace{ f(L) + \cdots + f(L)}_N.$$
\end{definition}

\begin{definition} The element $f$ is {\it strongly elliptic} on $L$ if there exist finite sets $M_i \subseteq \underbrace{ L \times \cdots \times L}_{n-1}$, $1 \leq i \leq n$, such that $Span(f(L))$ is a sum of additive subgroups 
$$f(a_1, \ldots, a_{i-1},L,a_{i+1}, \ldots, a_n),$$ where $(a_1, \ldots, a_{i-1}, a_{i+1}, \ldots,a_n) \in M_i$, $1 \leq i \leq n$.
\end{definition}

Clearly, if $f$ is strongly elliptic then it is elliptic.

To prove Theorem 3.9, in next section, the following two results are needed.

\begin{theorem} Let $L$ be a $Z/p^kZ$-Lie algebra generated by a finite set $X \subseteq L$. Let us assume that an arbitrary commutator in $X$ is ad-nilpotent, that is, the corresponding adjoint map is nilpotent (and elements in $X$ are considered as commutators of length 1).  Let $f(x_1, \ldots, x_n) \in Lie \langle X \rangle$ be a multilinear element in the free Lie $Z_p$-algebra such that $f \notin pLie \langle X \rangle$.  Then

i) $I = Span(f(L))$ is an ideal in $L$ and $|L:Span(f(L))| < \infty$,

ii) $f$ is strongly elliptic on $L$.

\end{theorem}

Now assume that $L= L(G)$ is the Lie algebra linked to the lower central series of a pro-$p$ group $G$ and consider 
$F$, the free  pro-$p$ group on the set of free generators $ \{x_1, x_2, \cdots \}$, and $L(F) = \oplus F_i/F_{i+1}$ the free Lie $Z_p$-algebra over the set of free generators $ \{ {\bar x_i} = x_iF_2 \; | \; i \geq 1\}$, as we mentioned above.

If $w = {\bar w}w' \in F$ is a multilinear word,   with
$${\bar w }= \prod_{\sigma \in S_n} [ \cdots[x_{\sigma(1)}, x_{\sigma(2)}], \cdots, x_{\sigma(n)}]^{k_{\sigma}}, k_{\sigma} \in Z_p, \; w' \in F_{n+1}, $$
then $f({\bar x_1}, \ldots, {\bar x_n}) \in L(F)$ is a multilinear element and $Span(f(L))$ is an ideal in $L$.

\begin{theorem}  With the above notation, if $|L: Span(f(L))| < \infty$ and $f$ is strongly elliptic on $L$, then 
$|G: \bar{\langle w(G) \rangle}| < \infty$  and $w$ is strongly elliptic on $G$.
\end{theorem}

In \cite{MM} it is proved that an arbitrary multilinear polynomial is strongly elliptic on the centerless Virasoro algebra {\it Vir} and on its subalgebras $Vir^{(k)}$, $k \geq -1$, where {\it Vir} has a $K$-basis $\{e_i \;| \;i \in Z \}$ with the multiplication $[e_i,e_j] = (i-j)e_{i+j}$, $K$ is a field of zero characteristic and $Vir^{(k)} = \sum_{i =k}^{\infty} K e_i$.

\medskip

\subsection{Finitely generated nil algebras}

In what follows, $K$ will denote a field of characteristic $\neq 2$.  Let $X =\{x_1, x_2, \ldots \}$ and let $K\langle X \rangle$ be the free associative $K$-algebra on the set of free generators $X$.

It is clear that for an arbitrary element $f \in K\langle X \rangle$ and an arbitrary $K$-algebra we can consider the notion of ellipticity and strong ellipticity of $f$ on any associative algebra $A$.

In the Lie algebra $K\langle X \rangle^{(-)}$ we can consider the Lie subalgebra $Lie\langle X \rangle$ generated by $X$.  We will refer to elements of $Lie\langle X \rangle$ as Lie polynomials.

The following results were proved in \cite {M2}.

\begin{lemma}  Let $\mathcal V $ be a  nontrivial variety of associative $K$-algebras.  There exists a nonzero Lie multilinear polynomial $h$ such that $h = 0$ in $\mathcal V$.
\end{lemma}

\begin{proposition}  For an arbitrary nonzero multilinear element $f \in K \langle X \rangle$, the space $Span f(K \langle X \rangle)$ contains a nonzero multilinear Lie polynomial.
\end{proposition}

\begin{theorem}  Let $A$ be a finitely generated nil $K$-algebra. Suppose that a nonzero multilinear element $f$ lies in $Lie \langle X \rangle$.  Then $f$ is strongly elliptic on $A$.
\end{theorem}

\begin{theorem} Let $A$ be a finitely generated nil $K$-algebra. An arbitrary nonzero multilinear element from $K \langle X \rangle$ is strongly elliptic on $A$.
\end{theorem}

\bigskip

\section{ Waring problems in associative algebras}

Let $K$ be a commutative unital ring and let
$f=f(X_1,\dots, X_m)$ be a (noncommutative) polynomial, i.e., an element of the free algebra
$F\langle X_1,X_2,\dots\rangle$. Let $A$ be an associative $K$-algebra. The set
$$f(A)= \{f(a_1,\dots,a_n)\,|\, a_i\in A\}$$
is called the {\it image of $f$} in $A$. By span$\,f(A)$ we denote the linear span of $f(A)$. 

Adapting some of the terminology introduced above, we make the following definition.

\begin{definition}
We  say that $f$ has {\it finite width in $A$} if there exists an $N\ge 1$ such that every element  in span$\,f(A)$ is a linear combination of $N$ elements from $f(A)$. The smallest such $N$ is  called the {\it width  of $f$ in $A$} and will be denoted by
$W_{f,A}$. If $f$ does not have finite width in $A$,  we write $W_{f,A}=\infty$. If every polynomial $f$ has finite width in $A$, then
 $$W_A:=\sup_f W_{f,A}$$ will be called  the {\it Waring constant of $A$}. 
\end{definition}

By  {\it Waring problems in $A$} we roughly mean   the study of the notions just introduced, along with their variations.

 \begin{remarks} 
   (a) It may be worthwhile  to  independently consider a more general situation  where "linear combination" in the above definition is replaced by "sum". Moreover, one may also consider the situation where $A$ is merely a subset of an algebra. In this way, one in particular covers the classical Waring problem where $f=X^n$, $K=\ZZ$, and $A=\NN\cup\{0\}$ 
    ($W_{f,A}$ is then the   function $g(n)$ mentioned in the introduction). Moreover, it covers  various  generalizations of the classical problem that ask whether  elements from certain rings can be presented as sums of a fixed number of $n$th powers. However, we will not discuss this topic here; we only refer to a few papers  \cite{KG, KVZ, LWo, PSS, Vas} where one can find further references. 
    
    (b) In some of the results mentioned below, the only coefficients in linear combinations that are needed are $\pm 1$. More precisely, some algebra elements will be expressed as sums of elements from the set $$f(A)-f(A) =\{s-t\,|\, s,t\in f(A)\}.$$ One should add, however, that having both $1$ and $-1$ at our disposal is  much more convenient than  allowing only $1$ as the   coefficient  (see  \cite{Vas}).
    \end{remarks}
		
		\medskip

\subsection{ Waring problems for general polynomials}
To tackle  Waring problems for arbitrary polynomials, it is natural to   first consider  matrix algebras. We begin by  stating a rather general result from \cite{B} which treats algebras of matrices  over an arbitrary algebra (over a field). Before that, we recall that a  polynomial $f$ is a polynomial
identity on an algebra $M$
if $f(M)=\{0\}$, and that $f$ is a central polynomial for $M$ if $f(M)$  lies in the center of $M$ but $f$ is not a polynomial identity of $M$. 

The following is \cite[Theorem 3.12]{B}.

\begin{theorem}\label{tt2}Let $K$ be an infinite field with {\rm char}$(K)\ne 2$, let $n\ge 2$, 
  let $f$ be a    
  polynomial  which is neither an identity nor a central polynomial of $M_n(K)$, let $k\ge 1$,
	and let $B$ be a unital $K$-algebra such that every element in $A=M_n(B)$ is a sum of $k$ commutators and a central element. 
	 Then every commutator in $A$ is 
	a sum of $1936k^2 + 22k$ elements from $f(A)-f(A)$. 
\end{theorem}

The proof combines  tools from PI theory, Lie theory of associative algebras, and results related to the presentation of commutators as sums of square-zero elements.

The number $1936k^2 + 22k$ surely is not optimal. The point of the theorem is that such a number exists and is independent of the size of the matrices, provided of course that so is $k$. This holds for various algebras $B$:
\begin{enumerate}
    \item[(a)] If $B$ is a commutative algebra,
    then every
    trace zero matrix in $M_n(B)$ is a sum of at most two commutators \cite{Me}. Therefore,
    assuming 
    that the characteristic of $K$ is  $0$, we see from
    $a=\left(a-\frac{{\rm tr}(a)}{n}1\right) + \frac{{\rm tr}(a)}{n}1$
    that the number $k$ from Theorem \ref{tt2} is (at most) $2$ for every $n$. Therefore,
    every commutator in $A=M_n(B)$ is a sum of 
    $1936\cdot 4 + 22\cdot 2 = 7788$ elements from $f(A)$.

  \item[(b)]
  Let $K$ be as in Theorem \ref{tt2}, let $V$ be an infinite-dimensional vector space over $K$, and let $A={\rm End}_K(V)$. Since $V^n$ is isomorphic to $V$ for every $n\in\NN$, $A$ is isomorphic to 
  $M_n(A)$. Given any nonconstant  polynomial $f$, we may choose $n$ such that $f$ is neither a polynomial identity nor a central polynomial of $M_n(K)$. Using the fact that every element in $A$ is a commutator \cite{Me}, it follows that every element is a sum of
  $1936\cdot 1 + 22\cdot 1 =1958$ elements from $f(A)$. This readily implies that $W_A\le 2\cdot 1958=3916$. Using a more refined method,  this $4$-digit number can be lowered to $14$, provided  that $K$ is algebraically closed  and has characteristic $0$  \cite[Corollary 3.11]{BV}. The exact value of the Waring constant $W_A$, however, is so far unknown.

 \item[(c)] Let $H$ be an infinite-dimensional Hilbert space and let $A=B(H)$, i.e., $A$ is the algebra of all bounded linear operators of $H$. As in (b), we have  $A\cong M_n(A)$ for every $n\in \NN$. Further,
  every element
 in $A$ that is not of the form $\lambda I + K$ with $\lambda\in\CC\setminus{\{0\}}$
and $K$ compact is a commutator \cite{BP}. In particular, compact operators are
commutators, so every element in $A$ is a sum of a commutator and a central
element $\lambda I$. That is, $k=1$.
  As every 
  element of $B(H)$ is a sum of two commutators \cite{Halmos}, it follows that the Waring constant $W_A$ exists and is bounded by a certain $4$-digit number.
  Again, more sophisticated methods yield better results. By
  \cite[Theorem 4.2]{BS1}, $W_A\le 8$. Moreover, \cite[Corollary 3.12]{BV} considers more general algebras $B(X)$ where $X$ is a Banach space such that $X\cong X\oplus X$. However, the exact value of $W_A$  is  unknown.
 \end{enumerate}

Let us now focus on the basic and most interesting case where our algebra is the ordinary matrix algebra $$M_n=M_n(\CC),\,\,\,n\ge 2.$$ 
The restriction to the field $\CC$ is made primarily for 
  simplicity of exposition. We ask the interested 
  reader to look at the original sources to check what conditions on the field are really required in the results   stated below.

Note first that 
$W_{f,M_n}=1$ if $f$ is either a polynomial identity or a central polynomial of $M_n$. We may therefore assume that $f$
is not such a polynomial. Then
span$\,f(M_n)$ can only be either the space sl$_n = {\rm sl_n}(\CC)$  of all trace zero matrices or the whole algebra $M_n$  \cite{BKlep}. Using this together with the fact that every trace zero matrix is a commutator \cite{Shoda}, it easily follows from Theorem \ref{tt2}  
that $W_{M_n} < \infty$.

Computing $W_{M_n}$ is another story. While Theorem \ref{tt2} gives a $4$-digit upper bound, it was shown already in \cite{B} that $W_{M_n}\le 9$.  A further development was made in 
\cite{BS1, BS2, BV}. We now review some of the results from these papers.

Let us start by 
pointing out that $W_{M_n}\ne 1$. In fact, $W_{X^2,M_n}\ne 1$. Indeed,
 a nilpotent matrix $a \in M_n$
of maximal nilindex 
is not a square of some matrix in $M_n$. This is because 
 $a=b^2$ 
implies $b^{2n}
 =a^n=0$ and hence
 $b^n
 = 0$, leading
to the contradiction that $a^{n-1} 
= b^{2n-2} = 0$.

It is less obvious that $W_{M_n}$ is not always equal to $2$.  
Roughly speaking,
$W_{f,M_n}$ is often  $2$, yet not for every $f$. Recall from the PI theory that $f$ is a $2$-central polynomial for $M_n$ if $f^2$ is a central polynomial, but $f$ is not. For instance, $[X_1,X_2]$ is a $2$-central polynomial for $M_2$. It is known that 
$2$-central polynomials also exist for some even $n >2$ \cite{KBMR2}.

The following is \cite[Theorem 1.2]{BS1}.

\begin{theorem}\label{pet}
 Suppose $n>2$ and there exists a $2$-central polynomial $f$ for $M_n$. Then 
$W_{f,M_n} \ge 3$, and hence $W_{M_n} \ge 3$.
\end{theorem}

On the other hand, \cite[Theorem 1.1]{BS1} implies the following.

\begin{theorem}\label{pet2}
$W_{M_n}\le 5$ for every $n\ge 2$.
\end{theorem}

Computing $W_{M_n}$ remains an open problem. 

 Theorem \ref{pet}
indicates that for a given polynomial $f$,
 pathological phenomena  may occur
for matrices of small sizes relative to degree of $f$. To avoid this, it is natural to determine $W_{f,M_n}$ for any large enough $n$. This  problem was first considered in \cite{BS2} and then continued in \cite{BV}. We will discuss only the latter paper which, neglecting some technical assumptions,  contains stronger results.

The following theorem is a special case of \cite[Theorem 3.2]{BV}.

\begin{theorem}
   \label{rnc} Let $f$ be a nonconstant polynomial. Then ${\rm sl}_n\subseteq f(M_n)-f(M_n)$ for all large enough $n$.
\end{theorem}

The proof uses tools from the theory of central simple algebras, algebraic geometry, free analysis, and classical number theory.
It is noteworthy to mention that  the theorem also holds in the more general situation where $f$ is a noncommutative rational function. 
The basic theorem in \cite{BV}, from which all others (including Theorem \ref{rnc}) are deduced, states that given a nonconstant noncommutative rational function $f$ in $m$ variables,
for any $n$ large enough there exist  matrices $A_1,\dots,A_m\in M_n$
 such that the matrix
$f(A_1,\dots,A_m)$ has $n$ distinct and nonzero eigenvalues.

Theorem \ref{rnc} gives a definitive solution to our problem for polynomials that are   sums of commutators (of some polynomials). This is because the image of such polynomials is  contained in sl$_n$. For polynomials that are not of that type one would wish to express any matrix in $M_n$ as a linear combination of two matrices from $f(M_n)$. We now state 
\cite[Theorem 3.5]{BV} which shows that 
this is "almost" true.

\begin{theorem}\label{tlike}
If a polynomial $f$ is not a sum of commutators, then, for $n$ large enough,
every nonscalar matrix in $M_n$ is a linear combination of two matrices from $f(M_n)$.  
\end{theorem}

We do not know  whether scalar matrices should really be excluded in this theorem.  What can be said as an immediate corollary is that scalar matrices are linear combinations of three matrices from $f(M_n)$ (so $W_{f,M_n}\le 3$ for all large enough $n$). However, perhaps $3$ can be replaced by $2$.

A natural and interesting problem is to extend some of the above results to finite-dimensional simple algebras. For finite-dimensional division algebras, everything is  open at present.
One of the questions that can be asked is the following:
Given a field $K$, does there exist a constant $C$ such that $W_A\le C$ for every finite-dimensional  division $K$-algebra $A$?
(We remark here that 
if  $A$ is any finite-dimensional $K$-algebra, then we trivially have
    $W_A \le \dim_K A$).

We close this subsection by mentioning the papers \cite{ChWa, PaPr} which consider  similar problems in the algebra of upper triangular matrices.

\medskip

\subsection{The L'vov-Kaplansky conjecture}

The problem to determine possible images of polynomials in matrix algebras $M_n(K)$ was supposedly initiated by I.\ Kaplansky.  In \cite{Dt}, I.\ V.\ L'vov proposed the following problem, which  is now  commonly named  after   him and Kaplansky.
\smallskip 

\noindent
{\it L'vov-Kaplansky conjecture:}   If $f$ is a multilinear polynomial,  $K$ is a field and $n\ge 2$, then $f(M_n(K))$ is a vector space.

\smallskip 

Using the notation introduced above, this can be reformulated as 
$$W_{f,M_n(K)}=1$$
whenever $f$ is multilinear.

As already mentioned above,
span\,$f(A)$ 
can only be equal to 
$\{0\}$ (if $f$ is a polynomial identity), the space $K$ of scalar matrices 
(if $f$ is a central polynomial),
the space of zero trace matrices   sl$_n(K)$, and the whole algebra $M_n(K)$. The L'vov-Kaplansky conjecture can therefore be  
 further reformulated as  $$f(M_n(K))\in\left\{\{0\}, K, {\rm sl}_n(K), M_n(K)\right\}$$
for every multilinear polynomial $f$.

Interest in this conjecture was revived around 15 years
ago by A.\ Kanel-Belov, S.\ Malev and L.\ Rowen \cite{KBMR} who proved the following theorem. 

\begin{theorem}
    If $K$ is a quadratically  closed field and $f$ is a multilinear polynomial,
    then $f(M_2(K))$ is a vector space.
\end{theorem}

In other words, the L'vov-Kaplansky 
conjecture is true for $n=2$ (and $K$ quadratically closed).
In spite of many efforts, the conjecture 
has remained open for  $n >2$.  We refer the reader to the survey paper \cite{KBMRY}, which provides a rather complete information on partial solutions.

Instead of considering matrices of small size, one can consider polynomials of small degree.
The following result was recently obtained by D.\ Vitas \cite{Vit3} (generalizing  the result of \cite{DK} which  covered matrices of even size).

\begin{theorem}
    If $K$ is an algebraically closed field of characteristic $0$ and $f$ is a multilinear polynomial of degree $3$ (or less), then 
    $f(M_n(K))$ is a vector space for every $n\ge 2$.
\end{theorem}

All in all, it is  fair 
to say that  we are far from
the solution of the L'vov-Kaplansky conjecture. On the other hand, many authors have  successfully studied different versions of  
this conjecture. 
Here is a list of some of the algebras $A$ having the property that 
$f(A)$ is a vector space for every multilinear polynomial $f$:
\begin{enumerate}
    \item[(a)] $A=\mathbb H$, the algebra of quaternions \cite{Mal}.
    \item[(b)] $A$ is the algebra of all strictly upper triangular matrices over a field $K$ \cite{Fag}.
    \item[(c)] $A$ is the algebra of all upper triangular matrices over an infinite field $K$
    \cite{GdM, LW}.

\item[(d)] $A$ is a unital algebra over a   field $K$ having a surjective inner derivation \cite{Vit} (this includes Weyl algebras and the algebra End$_K(V)$ where $V$ is an infinite-dimensional vector space).
\end{enumerate}

For other related results,  see  \cite{CF,  Fag2, FagKo, KBMPR, Vit2} and references therein.

\medskip

\subsection{The commutator width}
 We call the width of the polynomial $[X_1,X_2]=X_1X_2 - X_2X_1$ in the algebra $A$, i.e.,
$$\gamma_A:=W_{[X_1,X_2],A},$$
the {\it commutator width of $A$}.

As usual, we write $[A,A]$
for the linear span of all commutators in $A$. The condition that
$\gamma_A\le  n$ thus means that every element in $[A,A]$ is a sum of at most $n$ commutators, and if $\gamma_A= n$ then some elements in $[A,A]$ cannot be written as sums of $n-1$ commutators. If
for every $n$ there is an element in $[A,A]$
that is not a sum of $n$ commutators, we have
 $\gamma_A=\infty$. For example, it is easy to see that this holds if
 $A$ is the free algebra $
K\langle X_1,X_2,\dots\rangle$ (see \cite[Lemma 4.4]{BGS}).

In 1937, K.\ Shoda \cite{Shoda} proved that 
every  matrix
in sl$_n(K)$, $n\ge 2$, is a commutator, provided that $K$ is a field of characteristic $0$. Note that this is just another way of saying that $\gamma_{M_n(K)}=1$. Twenty years later,   A.\ Albert and B.\ Muckenhoupt \cite{AM} proved that this is true for any field $K$.  More recently, A.\ Stasinski \cite{S} showed that it is enough to assume that
 $K$ is a  principal ideal domain.  On the other hand, 
 $\gamma_{M_2(K)} = 1$
 does not hold for 
every commutative unital ring $K$ \cite{RR}, but we always have
 $\gamma_{M_n(K)} \le 2$ \cite{Me, RR}.

We summarize these results as follows.

\begin{theorem}
    Let  $K$ be a commutative unital ring. 
Then $\gamma_{M_n(K)}\le 2$
for every $n\ge 2$ and  there   exist  commutative unital rings
$K$ such that $\gamma_{M_2(K)} = 2$. However, if $K$ is a principal ideal domain (in particular, if $K$ is a field), then 
 $\gamma_{M_n(K)} = 1$ for every $n\ge 2$.    
\end{theorem}

What is the 
commutator width of a  simple algebra $A$? 
The Albert-Muckenhoupt
theorem answers this question for the basic case where $A=M_n(K)$, $K$ a field. It does not seem easy to generalize this to more general simple algebras. The following was proved by 
S.\ A.\ Amitsur and L.\ Rowen in 1994 \cite{AR}.

\begin{theorem}
    If $A$ is a finite-dimensional simple algebra, then $\gamma_A\le 2$.
\end{theorem}

    No example with $\gamma_A=2$ is known. It  is thus  open  whether or not the commutator width of a finite-dimensional simple algebra is actually $1$.  Results in \cite{AR} show that this is true in some special situations.

The commutator width of infinite-dimensional simple algebras is not always $1$. The following theorem combines two results.
The first one was proved by L.\ 
Makar-Limanov \cite{ML} 
and the second one by L.\ Robert \cite{Rob2}.

\begin{theorem}\label{412}
   There exists a division algebra $D$ (infinite-dimensional over its center) such 
that $\gamma_D=2$. Moreover, there exists a simple algebra $A$ such that $\gamma_A=\infty$.
\end{theorem}

Makar-Limanov's division algebra $D$ also satisfies 
$D=[D,D]$. Thus, every element in $D$ is a sum of at most two commutators, and there are elements in $D$ that are not commutators.

Robert's simple algebra $A$ is actually a $C^*$-algebra, so the proof is not only algebraic.
The result to which we are referring to is \cite[Theorem 1.4]{Rob2}. From its statement, however,  it  may not be immediately clear that
$\gamma_A=\infty$. Let us, therefore,  clarify. It is immediate that this theorem  implies that for every $n\in\NN$ there exists a self-adjoint element $a$ in $ [A,A]$ that cannot be written as $\sum_{i=1}^n [x_i^*,x_i]$ for some $x_i\in A$. 
Note that for any
self-adjoint elements
$s$ and $t$, we have
 $i[s,t]=[x^*,x]$ where $x=\frac{1}{\sqrt{2}}(s+it)$. Since every
 element in $A$ can be written as $u+iv$ with $u,v$ self-adjoint, and hence every commutator
 in $A$ can be written
 as $[x_1^*,x_1]+[x_2^*,x_2] + i([x_3^*,x_3]+[x_4^*,x_4])$ for some $x_i\in A$,
$\gamma_A=\infty$ follows.

Robert's theorem is interesting also from the Lie algebra theory viewpoint. Namely,  the $C^*$-algebra $A$  from this theorem has a tracial state and hence its center has trivial intersection with $[A,A]$, so the Lie algebra $[A,A]$ is  simple  (see \cite{Her}) and has infinite commutator width. This answers
\cite[Question 6]{DKR}.

There has actually been a great deal of work on sums of commutators in $C^*$-algebras. Let us mention only the fundamental result  by C.\ Pop \cite{Pop}, which can be in   our notation stated as follows: If $A$ is a unital $C^*$-algebra with no tracial states, then $[A,A]=A$ and $\gamma_A<\infty$.
 
\medskip

\subsection{The width of the product of commutators}

Another specific polynomial that has been systematically studied 
is $$h:=[X_1,X_2][X_3,X_4],$$ the product of two commutators. Let us write
$$\xi_A:=W_{h,A},$$
i.e., $\xi_A$ is the width of $[X_1,X_2][X_3,X_4]$ in $A$.

This study, however, began just recently.
In the seminal paper \cite{GT}, E.\ Gardella and H.\ Thiel proved 
the following  theorem.

\begin{theorem}\label{TGT}
    Let $A$ be a unital algebra over a commutative unital ring $K$. If $A$ is equal to its commutator ideal $A[A,A]A$,
    then $A={\rm span}\,h(A)$ and
    $\xi_A <\infty$.
\end{theorem}

This theorem in particular shows that
  $h$ has finite width in every unital $K$-algebra $A$ such that $A= {\rm span}\,h(A)$.  It is easy to see that the same is true for the polynomial $X^k$ provided that $K$ is  a field of  characteristic $0$ \cite{Vas}. One can now ask 
     what other polynomials $f$ have this property, that is, the property that $f$ has finite width in every unital $K$-algebra $A$ such that $A= {\rm span}\,f(A)$. We believe this problem is largely open.

Another problem that obviously stems from Theorem \ref{TGT} is to determine $\xi_A$ in algebras satisfying
$A=A[A,A]A$. We  briefly discuss it in the rest of this subsection.

Quite a while ago, J.\ D.\ Botha proved that
$\xi_{M_n(K)}=1$ for every field $K$ and every $n\ge 2$ \cite[Theorem 4.1]{Botha}.  More general matrix algebras were considered in \cite{GT}
and in the subsequent paper \cite{BGT}. We will not mention all results from these papers, but only
\cite[Theorem 5.4]{GT},
\cite[Theorem 2.2]{BGT},
and \cite[Theorem 4.4]{BGT}. We combine them in a unified statement.

\begin{theorem}\label{414}
    If $B$ is any unital algebra then
    $\xi_{M_n(B)}\le 2$ for every $n\ge 2$, and there exists a unital commutative algebra $B$ such that  $\xi_{M_n(B)}= 2$.
    If $D$ is a division algebra with infinite center, then  $\xi_{M_n(D)}= 1$ for every $n\ge 2$.
\end{theorem}

The question whether the assumption on the infinity of the center is necessary is  open. 

What are the possible values of $\xi_A$ for simple unital algebras $A$? This question was raised in \cite{GT} and addressed in 
the recent paper \cite{BJR}. We gather the main results of \cite{BJR} in one statement.

\begin{theorem}
    If $A$ is a finite-dimensional division algebra, then $\xi_A=1$. On the other hand, for every $n\in\NN$ there exists a simple unital algebra $A_n$ (infinite-dimensional over its center) such that $\xi_{A_n} > n$.
\end{theorem}

Algebras $A_n$ are $C^*$-algebras and the proof of the second statement of Theorem \ref{414}, just like the proof of the second statement of Theorem \ref{412},
rests on Villadsen's techniques from \cite{Vil}. The result that is proved in \cite{BJR}  actually  tells more than  just $\xi_{A_n} > n$.  It states that the unity $1$ of $A_n$ cannot be written as $\sum_{i=1}^n x_i[y_i,z_i]w_i$ for some $x_i,y_i,z_i,w_i\in A_n$.

E.\ Gardella and H.\ Thiel showed that $\xi_D\le 2$ for every division algebra $D$ \cite[Proposition 5.8]{GT}, and asked whether
$\xi_D$ is actually equal to $1$. The first statement of Theorem  \ref{414} gives a positive answer in the finite-dimensional case. In general, the question remains open.

\medskip

\subsection{Multiplicative Waring problems} 

So far, we have considered linear combinations of elements from the images of polynomials.
Therefore, a more accurate name for the problems that we considered would have been the {\it linear Waring problems}. In this last subsection, we 
consider the {\it multiplicative Waring problems}. By this we mean problems of expressing elements as 
products of elements from the  images of polynomials.
This topic was initiated in \cite{BV}. We will state two results from this paper. 

As above, we write $M_n$ for $M_n(\CC)$. Our first result is \cite[Theorem 4.6]{BV}.

    \begin{theorem}\label{lcm}
        Let $f$ and $g$ be nonconstant polynomials. For all $n$ large enough, every nonscalar invertible matrix in $M_n$ 
        is the  product of a matrix from $f(M_n)$ with a matrix from $g(M_n)$ (i.e., ${\rm GL}_n(\CC)\setminus\CC \subseteq f(M_n)\cdot g(M_n)$).
    \end{theorem}
    
Like in Theorem \ref{tlike}, we  do not know  whether the exclusion of scalar matrices is necessary.

The second result is \cite[Theorem 4.9]{BV}. Unlike Theorem \ref{lcm}, which considers only invertible matrices, it considers all matrices in $M_n$ (which is why the number of factors is considerably larger). 

\begin{theorem}
     If a nonzero polynomial $f$ has zero constant term,  then, for all $n$ large enough,
every matrix in $M_n$ is a product of twelve matrices from $f(M_n)$.   
\end{theorem}

The assumption that $f$
has zero constant term can be replaced by the milder assumption that
$f$ has a root in $\CC$. However, we do not know whether any assumption on $f$ (besides being nonconstant) is  necessary. Also, we do not know what is the minimal number that can replace twelve in this statement.

\end{document}